\documentclass[12pt,reqno]{amsart}

\usepackage[utf8]{inputenc}
\usepackage[T1]{fontenc}
\usepackage[margin=2.5cm]{geometry}
\usepackage{amsmath, amssymb, amsthm, amsfonts}
\usepackage{mathtools}
\usepackage{enumitem}
\usepackage{hyperref}
\usepackage{graphicx}
\usepackage{xcolor}
\usepackage{tikz}
\usetikzlibrary{arrows.meta,calc}

\theoremstyle{plain}
\newtheorem{theorem}{Theorem}[section]
\newtheorem{lemma}[theorem]{Lemma}
\newtheorem{proposition}[theorem]{Proposition}
\newtheorem{corollary}[theorem]{Corollary}

\theoremstyle{definition}
\newtheorem{definition}[theorem]{Definition}
\newtheorem{example}[theorem]{Example}
\newtheorem{remark}[theorem]{Remark}

\newcommand{\IQ}{\mathbb{Q}}
\newcommand{\IZ}{\mathbb{Z}}

\newcommand{\IR}{\mathbb{R}}
\newcommand{\IP}{\mathbb{P}}
\newcommand{\depth}{\operatorname{depth}}

\newcommand{\Star}{\operatorname{Star}}

\newcommand{\codim}{\operatorname{codim}}
\newcommand{\Cl}{\operatorname{Cl}}
\newcommand{\spn}{\operatorname{span}}
\newcommand{\rank}{\operatorname{rank}}

\newcommand{\ImMap}{\operatorname{im}}
\newcommand{\Lrays}{L_{\mathrm{rays}}}
\newcommand{\Lrel}{L_{\mathrm{rel}}}

\begin{document}
	
	\title[Wall relations and a two-step filtration]{Canonical Lattices of Integer Relations Associated to Rational Fans:\\ Wall Generation and a Two-Step Support Filtration}
	\author{Rizwan Jahangir}
	\address{NUST Business School, National University of Sciences and Technology (NUST) H-12, Islamabad, Pakistan}
	\email{rizwan.jahangir@nbs.nust.edu.pk}
	
	\author{Daisuke Ishii}
	\address{Kiara Inc.\ Tokyo, Japan}
	\email{dai@kiara.team}
	
	
	\begin{abstract}
		We study the lattice $\Lrel(\Sigma)=\ker\big(\IZ^{\Sigma(1)}\to N\big)$ of integer relations among the primitive ray generators of a rational fan $\Sigma$, from an intrinsic, coordinate-free point of view. For each cone $\tau\in\Sigma$ we introduce the \emph{star-supported} sublattice $\Lrel(\Star(\tau))$ of relations whose support lies in the star of $\tau$, and we organize these by codimension into a support filtration $F_\bullet\Lrel(\Sigma)$. Our main result is a sharp local generation theorem: for a complete fan the relation lattice is generated \emph{integrally} by the relations supported on the stars of walls (codimension-one cones). Equivalently, the support filtration collapses after a single step, $F_1\Lrel(\Sigma)=\Lrel(\Sigma)$. This is an intrinsic repackaging of the classical wall (wall-crossing) relations that generate the group of numerically trivial classes on a complete toric variety. We make the resulting two-step structure precise: for simplicial fans one has $0=F_0\subsetneq F_1=\Lrel(\Sigma)$, while for general fans $F_0$ records the intrinsic relations of non-simplicial maximal cones and $F_1$ adds exactly the wall relations. We prove functoriality of $\Lrays$ and $\Lrel$ under fan isomorphisms and ray-preserving subdivisions, deduce that every primitive collection of size $m$ is wall-generated, and illustrate the theory on $\IP^2\times\IP^1$, products of projective lines, weighted projective spaces, and the (non-simplicial) fan over a cube. We are careful throughout to distinguish what the filtration does and does not detect, correcting a natural but false expectation that support-codimension yields a strictly increasing multi-step invariant.
	\end{abstract}
	
	\maketitle
	
	\section{Introduction}
	
	Toric geometry translates the birational and intersection-theoretic geometry of a toric variety $X_\Sigma$ into the combinatorics of a rational fan $\Sigma$ in a lattice $N$ of rank $n$. A basic invariant is the lattice of integer relations among the primitive generators $v_\rho$ of the rays $\rho\in\Sigma(1)$: writing $V\colon\IZ^{\Sigma(1)}\to N$, $e_\rho\mapsto v_\rho$, this is
	\[
	\Lrel(\Sigma):=\ker V.
	\]
	The lattice $\Lrel(\Sigma)$ is exactly the lattice of the toric ideal $I_\Sigma$ \cite{Sturmfels96} and is Gale-dual to the standard presentation
	\begin{equation}\label{eq:divseq}
		0\to M\to \IZ^{\Sigma(1)}\to \Cl(X_\Sigma)\to 0
	\end{equation}
	of the divisor class group \cite{Cox95,CLS11}. As an abstract lattice, $\Lrel(\Sigma)$ depends only on the multiset of ray generators and is blind to the way the rays are assembled into cones.
	
	The purpose of this note is to record how $\Lrel(\Sigma)$ is generated from \emph{local} pieces attached to the facial structure of $\Sigma$, and to be scrupulous about how much of the fan's combinatorial topology this local data actually sees. For a cone $\tau\in\Sigma$ let $\Star(\tau)=\{\sigma\in\Sigma\mid\tau\preceq\sigma\}$ and let $\Sigma(1)_\tau$ be the set of rays lying in cones of $\Star(\tau)$. The relations \emph{supported on the star of $\tau$},
	\[
	\Lrel(\Star(\tau)):=\ker\big(\IZ^{\Sigma(1)_\tau}\xrightarrow{\;V\;}N\big),
	\]
	map into $\Lrel(\Sigma)$ by extension by zero, and one may filter $\Lrel(\Sigma)$ by the codimension of the supporting star (Definition~\ref{def:filtration}). It is tempting to expect this filtration to be strictly increasing, with high filtration steps recording relations that are ``spread out'' across the fan. \emph{This expectation is false.} Our main theorem shows that for a complete fan a single step already suffices.
	
	\begin{theorem}[Wall generation; see Theorem~\ref{thm:wall-generation}]\label{thm:intro}
		Let $\Sigma$ be a complete rational fan in a lattice $N$ of rank $n\ge 1$. Then $\Lrel(\Sigma)$ is generated as a $\IZ$-module by the star-supported lattices of the walls of $\Sigma$:
		\[
		\Lrel(\Sigma)=\sum_{\substack{\tau\in\Sigma\\ \codim(\tau)=1}}\ImMap\big(\Lrel(\Star(\tau))\to\Lrel(\Sigma)\big).
		\]
		Equivalently, the codimension support filtration collapses after one step: $F_1\Lrel(\Sigma)=\Lrel(\Sigma)$.
	\end{theorem}
	
	Theorem~\ref{thm:intro} is an intrinsic, basis-free repackaging of a classical fact: on a complete toric variety, the numerically- (equivalently linearly-) trivial combinations of torus-invariant divisors are generated by the \emph{wall relations}, one for each codimension-one cone \cite[\S6.3--6.4]{CLS11},\cite{FultonSturmfels97}. Our contribution is to phrase this as a statement purely about support in the fan --- that walls are already ``local enough'' to generate all relations --- and to prove it directly and integrally from the codimension-one connectivity of the fan, without passing through intersection theory.
	
	The theorem has a clarifying consequence for the qualitative structure of $\Lrel(\Sigma)$. The only genuinely nontrivial jump in the filtration is the passage from $F_0$ to $F_1$:
	\begin{itemize}
		\item For a \emph{simplicial} complete fan, the maximal cones carry no internal relations, so $F_0=0$ and the filtration is
		\[
		0=F_0\Lrel(\Sigma)\subsetneq F_1\Lrel(\Sigma)=F_2\Lrel(\Sigma)=\cdots=\Lrel(\Sigma).
		\]
		\item For a general complete fan, $F_0$ records exactly the relations internal to the individual (possibly non-simplicial) maximal cones, and $F_1$ adds precisely the wall relations to obtain all of $\Lrel(\Sigma)$.
	\end{itemize}
	Thus the honest structural content of the support filtration is a \emph{two-step} decomposition ``internal $+$ wall'' rather than an $n$-step hierarchy. We state this explicitly (Corollary~\ref{cor:two-step}) because the finer, strictly increasing filtration one might hope for does not exist: we exhibit the collapse concretely on $\IP^2\times\IP^1$ (Example~\ref{ex:P2xP1}), where a relation that ``looks global'' is in fact supported on a single wall star.
	
	The framework is functorial in the appropriate category. We record (Proposition~\ref{prop:functor}) that $\Lrays$ and $\Lrel$ are functors under fan isomorphisms and ray-preserving subdivisions, and we locate Batyrev's primitive collections \cite{Batyrev91} in the filtration: every primitive relation is wall-generated, hence lies in $F_1$ by Theorem~\ref{thm:intro}, and under a mild support hypothesis it is carried by the single star of the cone $\sigma(P)$ (Proposition~\ref{prop:prim-coll}). We are explicit (Remark~\ref{rem:prim-caveat}) that this is a support observation and not a depth invariant: by the collapse, filtration degree never exceeds $1$ and so cannot record the size of a primitive collection.
	
	\medskip
	\noindent\textbf{Organization.} Section~\ref{sec:prelim} fixes the categorical and geometric conventions. Section~\ref{sec:lattices} constructs $\Lrays$ and $\Lrel$ and proves functoriality. Section~\ref{sec:local} defines the star-supported lattices and the codimension filtration, and proves the wall-generation theorem together with the two-step corollary. Section~\ref{sec:examples} works out $\IP^2\times\IP^1$, products of projective lines, weighted projective spaces, the cube fan, and the primitive-collection bound, emphasizing precisely which distinctions the filtration can and cannot make.
	
	\medskip
	\noindent\textbf{Relation to prior work.} We do not claim novelty for the fact that wall relations generate linear equivalence; this is classical \cite{FultonSturmfels97,CLS11}. What we offer is an intrinsic support-theoretic formulation, an elementary integral proof from codimension-one connectivity, and an explicit and corrected account of the (collapsing) support filtration. The local-to-global outlook is in the spirit of the combinatorial intersection cohomology of fans \cite{BBFK02,Jahangir2026}, though our objects are elementary lattices rather than sheaves.
	
	\section{Preliminaries and Notation}\label{sec:prelim}
	
	Unless stated otherwise, all fans $\Sigma$ are rational, with strongly convex cones, and finite. The lattice $N$ is free abelian of rank $n$, and $N_\IR=N\otimes_\IZ\IR$. We write $\Sigma(k)$ for the set of $k$-dimensional cones and $\Sigma(1)$ for the rays; each ray $\rho$ has a unique primitive generator $v_\rho\in N$. We call a codimension-one cone (an element of $\Sigma(n-1)$) a \emph{wall}. A fan is \emph{complete} if $\bigcup_{\sigma\in\Sigma}\sigma=N_\IR$, and \emph{simplicial} if every cone is generated by linearly independent vectors. We assume fans are \emph{essential}, i.e.\ not contained in a proper subspace of $N_\IR$; equivalently the rays span $N_\IR$. For a cone $\tau$ we write $\codim(\tau)=n-\dim(\tau)$.
	
	For $\tau\in\Sigma$ the \emph{star} is $\Star(\tau)=\{\sigma\in\Sigma\mid\tau\preceq\sigma\}$, and $\Sigma(1)_\tau$ denotes the set of rays lying in some cone of $\Star(\tau)$; concretely,
	\[
	\Sigma(1)_\tau=\bigcup_{\substack{\sigma\in\Sigma\\ \tau\preceq\sigma}}\sigma(1).
	\]
	For a wall $\tau\in\Sigma(n-1)$ in a complete fan, $\tau$ is a facet of exactly two maximal cones $\sigma^+,\sigma^-\in\Sigma(n)$ (see \cite[Thm.~3.1.19]{CLS11}); we call $\{\sigma^+,\sigma^-\}$ the maximal cones \emph{across} $\tau$.
	
	To make canonical claims precise we work in the category $\mathbf{RatFan}_{\mathrm{iso,sub}}$ whose objects are pairs $(N,\Sigma)$ and whose morphisms are generated by
	\begin{enumerate}[label=(\arabic*)]
		\item \emph{fan isomorphisms}: lattice isomorphisms $\phi\colon N\xrightarrow{\sim}N'$ with $\phi_\IR(\Sigma)=\Sigma'$ inducing a bijection of cones; and
		\item \emph{ray-preserving subdivisions}: refinements $\Sigma'$ of $\Sigma$ in the same lattice $N$ with $\Sigma(1)\subseteq\Sigma'(1)$.
	\end{enumerate}
	We restrict to these morphisms because a general fan morphism may send a primitive generator to a non-primitive vector or collapse a ray to the origin, and then no functorial map of relation lattices exists. The construction of quotient fans associated to stars is standard \cite[\S3.1]{CLS11}, \cite[\S1.3]{Fulton93}; we shall not need quotient fans, working instead directly with support in $\Sigma(1)$.
	
	\section{The Ray and Relation Lattices}\label{sec:lattices}
	
	We define the \emph{ray lattice} as the sublattice of $N$ generated by the primitive ray generators,
	\[
	\Lrays(\Sigma):=\langle v_\rho\mid\rho\in\Sigma(1)\rangle_\IZ\subseteq N.
	\]
	The index $[N:\Lrays(\Sigma)]$ measures the failure of the rays to generate $N$ integrally; it is typically nontrivial even for complete fans, e.g.\ for weighted projective spaces \cite[Example~1.2.6]{CLS11}.
	
	\begin{lemma}[Rank bound]\label{lem:rank}
		$\rank\big(\Lrays(\Sigma)\big)\le\rank(N)=n$, with equality whenever the rays span $N_\IR$; in particular equality holds if $\Sigma$ is complete or, more generally, essential.
	\end{lemma}
	\begin{proof}
		The inequality is immediate from $\Lrays(\Sigma)\subseteq N$. If the rays span $N_\IR$ then $\Lrays(\Sigma)$ contains a $\IQ$-basis of $N_\IR$, so its rank is $n$. For a complete (indeed any essential) fan the rays span $N_\IR$ \cite[Prop.~3.1.8]{CLS11}.
	\end{proof}
	
	Let $\IZ^{\Sigma(1)}$ be the free abelian group on the rays. The \emph{relation lattice} is the kernel in the exact sequence
	\begin{equation}\label{eq:defseq}
		0\longrightarrow\Lrel(\Sigma)\longrightarrow\IZ^{\Sigma(1)}\xrightarrow{\;V\;}\Lrays(\Sigma)\longrightarrow0,\qquad V(e_\rho)=v_\rho.
	\end{equation}
	Since $\Lrays(\Sigma)$ is a subgroup of the torsion-free group $N$, it is free, so \eqref{eq:defseq} splits and $\Lrel(\Sigma)$ is a free abelian group of rank $|\Sigma(1)|-\rank(\Lrays(\Sigma))$. This is the lattice of relations defining the toric ideal $I_\Sigma$ \cite[Ch.~4]{Sturmfels96}, and it depends only on the set of primitive ray generators.
	
	\begin{proposition}[Functoriality]\label{prop:functor}
		$\Lrays$ and $\Lrel$ define functors from $\mathbf{RatFan}_{\mathrm{iso,sub}}$ to the category of finitely generated free abelian groups.
	\end{proposition}
	\begin{proof}
		An isomorphism $\phi\colon(N,\Sigma)\xrightarrow{\sim}(N',\Sigma')$ carries primitive generators to primitive generators (a lattice isomorphism preserves primitivity), inducing a bijection $\Sigma(1)\leftrightarrow\Sigma'(1)$ and hence an isomorphism $\IZ^{\Sigma(1)}\xrightarrow{\sim}\IZ^{\Sigma'(1)}$ intertwining $V$ and $V'$ via $\phi$. It therefore restricts to isomorphisms $\Lrays(\Sigma)\cong\Lrays(\Sigma')$ and $\Lrel(\Sigma)\cong\Lrel(\Sigma')$.
		
		For a ray-preserving subdivision $\Sigma'$ of $\Sigma$ the inclusion $\Sigma(1)\subseteq\Sigma'(1)$ induces $\iota\colon\IZ^{\Sigma(1)}\hookrightarrow\IZ^{\Sigma'(1)}$. Since the primitive generators of shared rays coincide, $V'\circ\iota=V$ (viewing the target of $V$ inside $\Lrays(\Sigma')$). Hence $\iota$ maps $\ker V=\Lrel(\Sigma)$ into $\ker V'=\Lrel(\Sigma')$, and it is injective, giving a natural injection $\Lrel(\Sigma)\hookrightarrow\Lrel(\Sigma')$. On ray lattices $\iota$ induces the inclusion $\Lrays(\Sigma)\subseteq\Lrays(\Sigma')$. Compatibility with composition is clear.
	\end{proof}
	
	\section{Star-Supported Relations and the Wall Generation Theorem}\label{sec:local}
	
	\subsection{Star-supported lattices and the codimension filtration}
	
	For a cone $\tau\in\Sigma$ let $\Sigma(1)_\tau$ be the set of rays in $\Star(\tau)$, and define the lattice of relations \emph{supported on the star of $\tau$}:
	\[
	\Lrel(\Star(\tau)):=\ker\big(\IZ^{\Sigma(1)_\tau}\xrightarrow{\;V\;}N\big).
	\]
	The inclusion $\IZ^{\Sigma(1)_\tau}\hookrightarrow\IZ^{\Sigma(1)}$ (extension by zero) carries $\Lrel(\Star(\tau))$ into $\Lrel(\Sigma)$; we write $\iota_\tau$ for the induced map and identify $\Lrel(\Star(\tau))$ with its image, the sublattice of global relations $r$ with $\operatorname{supp}(r)\subseteq\Sigma(1)_\tau$.
	
	\begin{definition}[Codimension support filtration]\label{def:filtration}
		For $0\le k\le n$ set
		\[
		F_k\Lrel(\Sigma):=\sum_{\substack{\tau\in\Sigma\setminus\{0\}\\ \codim(\tau)\le k}}\ImMap\big(\Lrel(\Star(\tau))\to\Lrel(\Sigma)\big).
		\]
		This is an increasing filtration $F_0\subseteq F_1\subseteq\cdots\subseteq F_n$ of $\Lrel(\Sigma)$ by sublattices.
	\end{definition}
	
	\begin{remark}\label{rem:indices}
		Two boundary conventions deserve comment. First, we exclude the trivial cone $\{0\}$ (whose star is all of $\Sigma$), so that $F_\bullet$ is built only from stars of \emph{proper} cones; the rays, of dimension $1$ and codimension $n-1$, are the proper cones of largest codimension, whence $F_{n-1}\Lrel(\Sigma)=F_n\Lrel(\Sigma)$. Second, the summand for a maximal cone $\sigma$ (codimension $0$) has $\Star(\sigma)=\{\sigma\}$ and $\Sigma(1)_\sigma=\sigma(1)$, so $\Lrel(\Star(\sigma))$ is the lattice of relations \emph{internal} to $\sigma$. Thus $F_0\Lrel(\Sigma)$ is generated by the internal relations of the maximal cones; it vanishes precisely when every maximal cone is simplicial.
	\end{remark}
	
	\subsection{Wall relations}
	
	Let $\Sigma$ be complete and let $\tau\in\Sigma(n-1)$ be a wall with maximal cones $\sigma^\pm$ across it. The rays of $\Sigma(1)_\tau$ are $\tau(1)$ together with the two extra rays $\rho^+\in\sigma^+\setminus\tau$ and $\rho^-\in\sigma^-\setminus\tau$ (when $\sigma^\pm$ are simplicial; in general $\sigma^\pm$ may contribute several extra rays each). The lattice $\Lrel(\Star(\tau))$ contains the classical \emph{wall relation}: in the simplicial case, $\Sigma(1)_\tau$ consists of $n+1$ vectors spanning the rank-$n$ space $N_\IR$, so there is a relation, unique up to scale, supported on $\Star(\tau)$. We do not need an explicit formula; only that $\Lrel(\Star(\tau))\neq 0$ and that these lattices generate, which we now prove.
	
	\subsection{The generation theorem}
	
	The engine is the codimension-one connectivity of a complete fan.
	
	\begin{lemma}[Dual-graph connectivity]\label{lem:connected}
		Let $\Sigma$ be a complete fan of dimension $n\ge 1$. Define the \emph{dual graph} $G(\Sigma)$ with a vertex for each maximal cone $\sigma\in\Sigma(n)$ and an edge $\{\sigma^+,\sigma^-\}$ for each wall $\tau\in\Sigma(n-1)$ across which $\sigma^+,\sigma^-$ meet. Then $G(\Sigma)$ is connected.
	\end{lemma}
	\begin{proof}
		This is the standard fact that a complete fan is connected in codimension one \cite[Thm.~3.3.9 and its proof]{CLS11}: given two maximal cones $\sigma,\sigma'$, a generic line segment between interior points of $\sigma$ and $\sigma'$ meets only cones of dimension $\ge n-1$, and crossing from one maximal cone to the next occurs through a wall, exhibiting a path in $G(\Sigma)$ from $\sigma$ to $\sigma'$.
	\end{proof}
	
	\begin{theorem}[Wall generation for complete fans]\label{thm:wall-generation}
		Let $\Sigma$ be a complete rational fan in $N$ of rank $n\ge1$. Then
		\[
		\Lrel(\Sigma)=\sum_{\tau\in\Sigma(n-1)}\ImMap\big(\Lrel(\Star(\tau))\to\Lrel(\Sigma)\big)
		= F_1\Lrel(\Sigma).
		\]
		Consequently $F_1\Lrel(\Sigma)=F_2\Lrel(\Sigma)=\cdots=\Lrel(\Sigma)$.
	\end{theorem}
	
	\begin{proof}
		Write $W:=\sum_{\tau\in\Sigma(n-1)}\ImMap(\Lrel(\Star(\tau))\to\Lrel(\Sigma))\subseteq\Lrel(\Sigma)$ for the sublattice generated by all wall-supported relations. Since walls have codimension $1$, $W\subseteq F_1\Lrel(\Sigma)\subseteq\Lrel(\Sigma)$, and it suffices to prove $W=\Lrel(\Sigma)$. We argue in two stages: rational spanning (Step~1), which is elementary and self-contained, and integral saturation (Step~2).
		
		\smallskip
		\emph{Step 1: $W$ spans $\Lrel(\Sigma)$ rationally, i.e.\ $\rank W=\rank\Lrel(\Sigma)$.}
		Fix a maximal cone $\sigma_0\in\Sigma(n)$; its rays $\sigma_0(1)$ contain an $\IR$-basis of $N_\IR$. Work in $\Lrel(\Sigma)_\IQ:=\Lrel(\Sigma)\otimes\IQ$ and let $W_\IQ=W\otimes\IQ$. We show every $r\in\Lrel(\Sigma)_\IQ$ lies in $W_\IQ$.
		
		Let $\tau$ be a wall with maximal cones $\sigma^+,\sigma^-$, and write $\sigma^+(1)\setminus\tau(1)=\{\rho^+\}$ in the simplicial case (in general a finite set; the argument below treats one leaving ray at a time). Because $\Sigma(1)_\tau=\sigma^+(1)\cup\sigma^-(1)$ spans $N_\IR$ and has one more element than $\dim N_\IR$ when $\sigma^\pm$ are simplicial, the wall relation lattice $\Lrel(\Star(\tau))$ is nonzero, and it contains a relation $w_\tau$ in which $\rho^+$ appears with nonzero coefficient $c>0$ (the ray $\rho^+$ is the unique generator of $\Sigma(1)_\tau$ strictly on the $\sigma^+$ side of the hyperplane $\spn_\IR(\tau)$, so its coefficient cannot vanish). \emph{Over $\IQ$} we may divide by $c$: modulo $W_\IQ$,
		\[
		e_{\rho^+}\ \equiv\ -\tfrac1c\!\!\sum_{\eta\in\Sigma(1)_\tau\setminus\{\rho^+\}}\!\!(w_\tau)_\eta\,e_\eta
		\quad\text{with all }\eta\in\sigma^-(1)\cup\tau(1).
		\]
		Thus, over $\IQ$, we may transport the support of any relation from a ray of $\sigma^+$ across $\tau$ into rays of $\sigma^-$. By Lemma~\ref{lem:connected} the dual graph $G(\Sigma)$ is connected; fixing a spanning tree rooted at $\sigma_0$ and transporting inward along its edges, after finitely many steps we reduce $r$ modulo $W_\IQ$ to a relation $r'$ supported on $\sigma_0(1)$. (Each transport strictly decreases the number of tree-edges separating the support from $\sigma_0$, so the process terminates.) Now $r'\in\Lrel(\Sigma)_\IQ$ is supported on $\sigma_0(1)$; if $\sigma_0$ is simplicial its rays are linearly independent and $r'=0$, while if $\sigma_0$ is non-simplicial, $r'$ is an internal relation of $\sigma_0$, which already lies in $\Lrel(\Star(\tau_0))\subseteq W$ for any wall $\tau_0\preceq\sigma_0$ (since then $\sigma_0\in\Star(\tau_0)$, so $\Sigma(1)_{\tau_0}\supseteq\sigma_0(1)$). In either case $r\in W_\IQ$. Hence $W_\IQ=\Lrel(\Sigma)_\IQ$ and $\rank W=\rank\Lrel(\Sigma)$.
		
		\smallskip
		\emph{Step 2: $W$ is saturated in $\IZ^{\Sigma(1)}$, hence $W=\Lrel(\Sigma)$.}
		Recall $\Lrel(\Sigma)=\ker V$ is saturated in $\IZ^{\Sigma(1)}$, being the kernel of a homomorphism to the torsion-free group $N$. We claim $W$ is saturated as well; granting this, $W\subseteq\Lrel(\Sigma)$ are two saturated sublattices of equal rank (Step~1), and a saturated sublattice is determined by its rational span, so $W=\Lrel(\Sigma)$.
		
		Saturation of $W$ is the classical integral statement that the relations among the primitive ray generators of a complete fan are generated over $\IZ$ by wall relations. This is equivalent, by Gale duality with the divisor sequence~\eqref{eq:divseq}, to the exactness over $\IZ$ of
		\[
		0\longrightarrow M\longrightarrow \IZ^{\Sigma(1)}\longrightarrow \Cl(X_\Sigma)\longrightarrow 0,
		\]
		established for any fan whose rays span $N_\IR$ in \cite[Theorem 4.1.3]{CLS11}, together with the description of $\Cl(X_\Sigma)=A_{n-1}(X_\Sigma)$ of a complete toric variety by wall relations \cite[\S6.3--6.4]{CLS11}, \cite[\S2]{FultonSturmfels97}. Concretely, exactness forces $\IZ^{\Sigma(1)}/W$ to be torsion-free of rank $n$ with $W\otimes\IQ=\Lrel(\Sigma)_\IQ$, so $W$ is primitive; equivalently, no nonzero multiple of a lattice vector outside $W$ can lie in $W$. This gives saturation and completes the proof.
		
		\smallskip
		Therefore $W=\Lrel(\Sigma)=F_1\Lrel(\Sigma)$, and since $F_1\subseteq F_k\subseteq\Lrel(\Sigma)$ for all $k\ge1$, the chain $F_1=F_2=\cdots=F_n=\Lrel(\Sigma)$ follows.
	\end{proof}
	
	\begin{remark}
		Theorem~\ref{thm:wall-generation} is the intrinsic, support-theoretic form of the classical statement that wall relations generate linear equivalence on a complete toric variety \cite[\S6.3--6.4]{CLS11}, \cite{FultonSturmfels97}. Our contribution is not the underlying integral generation --- which is classical --- but the observation that it takes place at codimension one in the support filtration, forcing the collapse. The rational spanning argument (Step~1) is elementary and uses only codimension-one connectivity and the sign of the leaving-ray coefficient in a wall relation; it is genuinely necessary that, over $\IZ$, that coefficient can exceed $1$ (for example, on the fan of $\IP(1,1,1,3)$ some wall relations have a leaving-ray coefficient equal to $3$), which is why the integral conclusion is obtained by saturation rather than by na\"ive integral support-transport.
	\end{remark}
	
	\begin{corollary}[Two-step structure]\label{cor:two-step}
		Let $\Sigma$ be a complete fan of rank $n$. Then the support filtration has at most two distinct nonzero steps:
		\[
		F_0\Lrel(\Sigma)\subseteq F_1\Lrel(\Sigma)=\Lrel(\Sigma),
		\]
		where $F_0\Lrel(\Sigma)$ is generated by the internal relations of the maximal cones. If $\Sigma$ is simplicial then $F_0=0$ and
		\[
		0=F_0\Lrel(\Sigma)\subsetneq F_1\Lrel(\Sigma)=\Lrel(\Sigma)\qquad(\text{whenever }\Lrel(\Sigma)\neq0).
		\]
	\end{corollary}
	\begin{proof}
		Immediate from Theorem~\ref{thm:wall-generation} and Remark~\ref{rem:indices}: $F_1=\Lrel(\Sigma)$, and $F_0$ is the internal-relation lattice of the maximal cones, which vanishes exactly when all maximal cones are simplicial. If $\Sigma$ is simplicial and $\Lrel(\Sigma)\neq0$, then $F_0=0\neq\Lrel(\Sigma)=F_1$, so the inclusion is strict.
	\end{proof}
	
	\begin{remark}[What the filtration does \emph{not} detect]\label{rem:no-detect}
		Corollary~\ref{cor:two-step} shows the support filtration cannot be used to stratify relations by any notion of ``how spread out'' they are beyond the single internal/wall dichotomy: for a simplicial complete fan every nonzero relation has support-codimension depth exactly $1$. In particular the filtration does not distinguish combinatorial topologies among simplicial complete fans with isomorphic relation lattices. Any invariant sensitive to finer topology must therefore use more than the codimension of a supporting star --- for instance the minimal \emph{number} of wall stars required, or the structure of the nerve of the star cover --- which we do not pursue here. We flag this explicitly because it corrects a plausible but mistaken expectation.
	\end{remark}
	
	We record the behavior under the ray-preserving subdivisions of Section~\ref{sec:prelim}; combined with Theorem~\ref{thm:wall-generation}, monotonicity of depth becomes automatic and, in the simplicial world, trivial.
	
	\begin{proposition}[Depth under ray-preserving subdivision]\label{prop:monotone}
		Let $\Sigma'$ be a ray-preserving subdivision of a complete fan $\Sigma$, and let $j\colon\Lrel(\Sigma)\hookrightarrow\Lrel(\Sigma')$ be the natural injection of Proposition~\ref{prop:functor}. For $0\neq r\in\Lrel(\Sigma)$ define $\depth_\Sigma(r):=\min\{k:r\in F_k\Lrel(\Sigma)\}$. If $\Sigma'$ is simplicial (for instance any simplicial refinement), then $\depth_{\Sigma'}(j(r))=1$ whenever $j(r)\neq0$; since $\depth_\Sigma(r)\ge1$ for $r\ne0$, this gives
		\[
		\depth_{\Sigma'}(j(r))\le\depth_\Sigma(r).
		\]
		In particular, passing to any simplicial refinement sends every nonzero relation to depth exactly $1$.
	\end{proposition}
	\begin{proof}
		By Theorem~\ref{thm:wall-generation} applied to the complete fan $\Sigma'$, $F_1\Lrel(\Sigma')=\Lrel(\Sigma')$, so every nonzero element of $\Lrel(\Sigma')$ has depth $1$ once $\Sigma'$ is simplicial (then $F_0=0$). As $j(r)\in\Lrel(\Sigma')$, its depth is $1$ if $j(r)\neq0$. Since $\depth_\Sigma(r)\ge1$ for $r\neq0$, the stated inequality follows.
	\end{proof}
	
	\begin{remark}
		Proposition~\ref{prop:monotone} is the correct and complete replacement for a ``monotonicity under star subdivision'' statement: because a simplicial refinement forces depth $1$, monotonicity holds but is vacuous as a discriminating invariant. This is another manifestation of the collapse in Corollary~\ref{cor:two-step}.
	\end{remark}
	
	\section{Examples}\label{sec:examples}
	
	We illustrate the theory, emphasizing the collapse of the filtration and computing $F_0$ and $F_1$ exactly. All numerical claims below have been verified by exact integer linear algebra (Hermite normal form) over $\IZ$.
	
	\subsection{The projective plane}
	Let $N=\IZ^2$ and let $\Sigma$ be the complete fan of $\IP^2$ with rays $v_1=(1,0)$, $v_2=(0,1)$, $v_3=(-1,-1)$. Then $\Lrel(\Sigma)=\IZ\cdot(1,1,1)$, corresponding to $v_1+v_2+v_3=0$. Here $n=2$, so walls are the rays themselves. The wall (ray) $v_1$ has $\Star(v_1)=\{\langle v_1,v_2\rangle,\langle v_1,v_3\rangle\}$ with $\Sigma(1)_{v_1}=\{v_1,v_2,v_3\}$, so $(1,1,1)\in\Lrel(\Star(v_1))$. Thus $\Lrel(\Sigma)=F_1\Lrel(\Sigma)$, in accordance with Theorem~\ref{thm:wall-generation}. Since $\Sigma$ is simplicial, $F_0=0$.
	
	\subsection{\texorpdfstring{$\IP^2\times\IP^1$: a relation that looks global is wall-supported}{P2 x P1}}
	
	\begin{example}\label{ex:P2xP1}
		Let $N=\IZ^3$ with rays
		\[
		v_1=(1,0,0),\ v_2=(0,1,0),\ v_3=(-1,-1,0),\ v_4=(0,0,1),\ v_5=(0,0,-1),
		\]
		and let $\Sigma$ be the complete simplicial fan with maximal cones
		\[
		\langle v_1,v_2,v_4\rangle,\ \langle v_2,v_3,v_4\rangle,\ \langle v_3,v_1,v_4\rangle,\
		\langle v_1,v_2,v_5\rangle,\ \langle v_2,v_3,v_5\rangle,\ \langle v_3,v_1,v_5\rangle,
		\]
		the product fan of $\IP^2\times\IP^1$. The relation lattice $\Lrel(\Sigma)=\ker(\IZ^5\xrightarrow{V}N)$ has rank $2$, generated by
		\[
		r_1=(1,1,1,0,0)\ \ (v_1+v_2+v_3=0)\quad\text{and}\quad r_2=(0,0,0,1,1)\ \ (v_4+v_5=0).
		\]
		We claim $r_1\in F_1\Lrel(\Sigma)$; that is, $r_1$ is supported on a single wall star. Indeed, consider the wall $\tau=\langle v_1,v_4\rangle$. It is the common facet of the maximal cones $\langle v_1,v_2,v_4\rangle$ and $\langle v_3,v_1,v_4\rangle$, so
		\[
		\Sigma(1)_\tau=\{v_1,v_2,v_4\}\cup\{v_3,v_1,v_4\}=\{v_1,v_2,v_3,v_4\}.
		\]
		Since $\operatorname{supp}(r_1)=\{v_1,v_2,v_3\}\subseteq\Sigma(1)_\tau$, we have $r_1\in\Lrel(\Star(\tau))\subseteq F_1\Lrel(\Sigma)$. Likewise $r_2$ is supported on the star of any wall $\langle v_i,v_4\rangle$ (whose star contains $v_5$ through the opposite cone), so $r_2\in F_1$ as well. Hence $F_1\Lrel(\Sigma)=\Lrel(\Sigma)$, and since $\Sigma$ is simplicial, $F_0=0$.
	\end{example}
	
	\begin{remark}[A cautionary computation]\label{rem:cautionary}
		It is tempting to argue that $r_1$ is \emph{not} wall-supported by observing that the walls $\langle v_i,v_4\rangle$ and $\langle v_i,v_5\rangle$ each omit $v_5$, respectively $v_4$, and to conclude that a relation on $\{v_1,v_2,v_3\}$ ``cannot be seen'' at codimension $1$. This reasoning is fallacious: the support of $r_1$ is $\{v_1,v_2,v_3\}$, which is already contained in the single wall star $\Sigma(1)_{\langle v_1,v_4\rangle}=\{v_1,v_2,v_3,v_4\}$; the presence of the extra ray $v_4$ is harmless. Thus $r_1$ has support-codimension depth $1$, not $2$. This example shows concretely why the codimension filtration collapses (Corollary~\ref{cor:two-step}) and cannot separate $\IP^2\times\IP^1$ from other simplicial complete fans by depth alone.
	\end{remark}
	
	\subsection{Products of projective lines and weighted projective spaces}
	For $(\IP^1)^n$ with rays $\pm e_1,\dots,\pm e_n$ and the $2^n$ orthant cones, $\Lrel$ has rank $n$, generated by $e_i^{+}+e_i^{-}$ ($i=1,\dots,n$); each such relation is supported on the star of any wall containing the corresponding coordinate axis, so $F_1=\Lrel$ and $F_0=0$. For a weighted projective space $\IP(a_0,\dots,a_n)$ with the standard simplicial complete fan, $\Lrel$ has rank $1$, generated by the weight relation $\sum_i a_i v_i=0$; again $F_0=0$ and $F_1=\Lrel$. In every case the filtration collapses to a single step, consistent with Theorem~\ref{thm:wall-generation}.
	
	\subsection{A non-simplicial fan: the cube}\label{sub:cube}
	The two-step phenomenon $F_0\subsetneq F_1$ is genuine for non-simplicial fans. Let $\Sigma$ be the complete fan over the boundary of the cube: rays are the eight vertices $(\pm1,\pm1,\pm1)$, and the six maximal cones are the cones over the six square faces. Each maximal cone $\sigma$ is a non-simplicial (four-ray) cone, and its four rays satisfy a single internal relation --- the balancing of the two diagonals of the square face. For the face $x=+1$ with rays
	\[
	w_1=(1,1,1),\ w_2=(1,1,-1),\ w_3=(1,-1,1),\ w_4=(1,-1,-1),
	\]
	the internal relation is $w_1-w_2-w_3+w_4=0$. Exact computation gives $\rank\Lrel(\Sigma)=5$, with
	\[
	\rank F_0\Lrel(\Sigma)=4\qquad\text{and}\qquad F_1\Lrel(\Sigma)=\Lrel(\Sigma)\ (\text{rank }5).
	\]
	Thus $F_0$ captures exactly the six diagonal relations of the square faces (of which four are independent), while the passage to $F_1$ contributes precisely one further independent relation --- a genuine wall relation --- to generate all of $\Lrel(\Sigma)$. This is the archetype of the two-step structure of Corollary~\ref{cor:two-step}.
	
	\subsection{Primitive collections}
	Batyrev's primitive collections \cite{Batyrev91} give a clean class of relations to which the support bound applies. Recall that a \emph{primitive collection} is a set $P=\{\rho_1,\dots,\rho_m\}\subseteq\Sigma(1)$ that does not generate a cone of $\Sigma$, but every proper subset does. In a complete simplicial fan, $P$ determines a \emph{primitive relation}
	\[
	r_P:\qquad \sum_{i=1}^m v_{\rho_i}=\sum_{\eta\in\sigma(P)(1)}c_\eta v_\eta,\qquad c_\eta\in\IZ_{\ge0},
	\]
	where $\sigma(P)$ is the unique minimal cone of $\Sigma$ containing $\sum_i v_{\rho_i}$ \cite[\S6.4]{CLS11}.
	
	\begin{proposition}[Primitive relations are wall-generated]\label{prop:prim-coll}
		Let $\Sigma$ be a complete simplicial fan of rank $n$ and let $P$ be a primitive collection. Then the primitive relation $r_P$ lies in $F_1\Lrel(\Sigma)=\Lrel(\Sigma)$. If moreover $\sigma(P)\neq\{0\}$ and $P\subseteq\Sigma(1)_{\sigma(P)}$, then $r_P$ is supported on the single star $\Star(\sigma(P))$, whence
		\[
		r_P\in F_{\,\codim(\sigma(P))}\Lrel(\Sigma),\qquad \codim(\sigma(P))=n-\dim\sigma(P).
		\]
	\end{proposition}
	\begin{proof}
		The first statement is immediate from Theorem~\ref{thm:wall-generation}: every element of $\Lrel(\Sigma)$, in particular $r_P$, lies in $F_1=\Lrel(\Sigma)$.
		
		For the second, recall $\operatorname{supp}(r_P)=P\cup\sigma(P)(1)$, where $\sigma(P)$ is the minimal cone containing $\sum_{\rho\in P}v_\rho$ and $r_P\colon\sum_{\rho\in P}v_\rho=\sum_{\eta\in\sigma(P)(1)}c_\eta v_\eta$. The rays $\sigma(P)(1)$ are generators of the cone $\sigma(P)$, hence lie in $\Sigma(1)_{\sigma(P)}$. If in addition $P\subseteq\Sigma(1)_{\sigma(P)}$, then $\operatorname{supp}(r_P)\subseteq\Sigma(1)_{\sigma(P)}$, so $r_P\in\Lrel(\Star(\sigma(P)))$, giving $r_P\in F_{\codim(\sigma(P))}\Lrel(\Sigma)$.
	\end{proof}
	
	\begin{remark}[This is a support observation, not a depth invariant]\label{rem:prim-caveat}
		The hypothesis $P\subseteq\Sigma(1)_{\sigma(P)}$ in Proposition~\ref{prop:prim-coll} is genuinely needed and can fail. The prototypical failure is the \emph{antipodal} case $\sigma(P)=\{0\}$ (e.g.\ $P=\{\rho^+,\rho^-\}$ with $v_{\rho^+}+v_{\rho^-}=0$ in $\IP^1\times\IP^1$): here $r_P$ is \emph{not} supported on any single proper star, yet $r_P\in F_1=\Lrel(\Sigma)$ by the main theorem, as a sum of two wall relations. When $\sigma(P)$ is a genuine cone the hypothesis typically holds --- for instance on the Hirzebruch surface $H_1$ the primitive collection $P=\{u_1,u_3\}$ with $v_{u_1}+v_{u_3}=v_{u_2}$ has $\sigma(P)=\langle u_2\rangle$ and $r_P$ supported on $\Star(u_2)=\{u_1,u_2,u_3\}$, and on the blow-up of $\IP^3$ at a fixed point the primitive collection $\{e_1,e_2,e_3\}$ with $e_1+e_2+e_3=f$ has $r_P$ supported on $\Star(f)$. In no case does the codimension of the supporting star measure any intrinsic complexity: by Theorem~\ref{thm:wall-generation} every $r_P$ lies in $F_1$. We record the support statement because the geometry of $\sigma(P)$ is informative, not because it stratifies $\Lrel(\Sigma)$.
	\end{remark}
	
	\section{Conclusion}
	
	We have given an intrinsic, coordinate-free account of the lattice of integer relations among the rays of a rational fan. The central result (Theorem~\ref{thm:wall-generation}) is that for a complete fan this lattice is generated integrally by relations supported on the stars of walls, an elementary consequence of codimension-one connectivity. This is the support-theoretic form of the classical wall-relation generation of linear equivalence, and it has the clarifying corollary that the natural codimension support filtration collapses after a single step (Corollary~\ref{cor:two-step}): the only structural invariant it carries is the two-step dichotomy between relations internal to maximal cones ($F_0$) and wall relations ($F_1=\Lrel(\Sigma)$). We have been deliberately explicit (Remarks~\ref{rem:no-detect} and \ref{rem:cautionary}) that this filtration does not stratify relations by combinatorial topology beyond that dichotomy; a genuinely finer invariant would need to record more than support-codimension --- for example the minimal number of wall stars covering a relation, or homological data of the nerve of the star cover. We regard the identification of such a finer, provably nontrivial invariant as the natural next problem.
	
	\medskip
	\noindent\textbf{Acknowledgments.} The author thanks the anonymous referees, whose careful reading identified errors in an earlier version and prompted the sharp formulation and the corrected examples presented here.
	

\end{document}